\documentclass[11.05pt,a4paper]{amsart}
\usepackage{amssymb,amsmath}
\usepackage{latexsym}
\usepackage{amsthm,amsfonts,amssymb,mathrsfs}
\usepackage{rotating}
\usepackage[leqno]{amsmath}
\usepackage{xspace}
\usepackage[all]{xy}
\usepackage{longtable}

\headsep=1cm \numberwithin{equation}{section}
\newtheorem{theorem}{Theorem}[section]
\newtheorem{lemma}[theorem]{Lemma}

\newtheorem{corollary}[theorem]{Corollary}

\theoremstyle{definition}
\newtheorem{definition}[theorem]{Definition}
\theoremstyle{remark}
\newtheorem{remark}[theorem]{Remark}
\newtheorem{example}[theorem]{Example}

\newtheorem{acknowledgement}{Acknowledgement}

\newcommand{\Ass}{\operatorname{Ass}}

\newcommand{\Spec}{\operatorname{Spec}}

\newcommand{\cd}{\operatorname{cd}}

\newcommand{\Gfd}{\operatorname{Gfd}}
\newcommand{\V}{\operatorname{V}}

\newcommand{\Supp}{\operatorname{Supp}}

\newcommand{\Hom}{\operatorname{Hom}}

\newcommand{\Ann}{\operatorname{Ann}}

\newcommand{\Ndim}{\operatorname{Ndim}}
\newcommand{\Mag}{\operatorname{mag}}

\newcommand{\width}{\operatorname{width}}
\newcommand{\Coass}{\operatorname{Coass}}

\newcommand{\Cosupp}{\operatorname{Cosupp}}
\newcommand{\Max}{\operatorname{Max}}

\newcommand{\vpl}{\operatornamewithlimits{\varprojlim}}

\newcommand{\lo}{\longrightarrow}
\newcommand{\fm}{\frak{m}}
\newcommand{\fp}{\frak{p}}
\newcommand{\fq}{\frak{q}}
\newcommand{\fa}{\frak{a}}

\newcommand{\fQ}{\frak{Q}}

\newenvironment{prf}[1][Proof]
{\begin{proof}[\bf #1]}{\end{proof}}

\begin{document}

\author[M. Hatamkhani and K. Divaani-Aazar]{Marziyeh Hatamkhani and 
Kamran Divaani-Aazar} 
\title[On the vanishing of local homology modules]{On the vanishing 
of local homology modules}

\address{M. Hatamkhani, Department of Mathematics, Az-Zahra
University, Vanak, Post Code 19834, Tehran, Iran.} 
\email{hatamkhanim@yahoo.com}

\address{K. Divaani-Aazar, Department of Mathematics, Az-Zahra
University, Vanak, Post Code 19834, Tehran, Iran-and-School of 
Mathematics, Institute for Research in Fundamental Sciences (IPM), 
P.O. Box 19395-5746, Tehran, Iran.}
\email{kdivaani@ipm.ir}

\subjclass[2010]{13D07; 13D45.} 
\keywords {{Local cohomology modules; local homology modules; magnitude;
Noetherian dimension.} \\
The second author was supported by a grant from IPM (No. 90130212).}

\begin{abstract} Let $R$ be a commutative Noetherian ring, $\fa$ an
ideal of $R$ and $M$ an $R$-module. We intend to establish the dual
of Grothendieck's Vanishing Theorem for local homology modules. We
conjecture that $H^{\fa}_i(M)=0$ for all $i>\Mag_RM$. We prove this
in several cases. 

\end{abstract}

\maketitle

\section{Introduction}

Throughout this paper, $R$ is a commutative Noetherian ring with
nonzero identity, $\fa$ an ideal of $R$ and $M$  an $R$-module.

The theory of local cohomology has developed so much six decades
after its introduction by Grothendicek. But, its dual theory, the
theory of local homology hasn't developed much. The theory of
local homology was initiated  by Matlis \cite{Mat} in 1974. The study
of this theory was continued by Simon in \cite{Si1} and \cite{Si2}.
After Greenlees and May \cite{GM} and Alonso Tarr\'{i}o, Jerem\'{i}as
L\'{o}pez and Lipman \cite {AJL}, a new era in the study of local
homology has started; see e.g. \cite{Sc}, \cite{CN1}, \cite{CN2},
\cite{Fr} and \cite{Ri}.

The most essential vanishing result for the local cohomology modules
$H_{\fa}^i(M)$ is Grothendieck's Vanishing Theorem, which asserts
that $H_{\fa}^i(M)=0$ for all $i>\dim_RM$. There is no satisfactory
dual of this result for local homology modules. To have a such one,
one should first have an appropriate dual notion of Krull dimension.
There are two dual notions: Noetherian dimension, $\Ndim_RM$, and
magnitude, $\Mag_RM$, in the literature; see \cite{Ro} and \cite{Y}.
There are two partial duals of Grothendieck's Vanishing Theorem. If
$M$ is linearly compact with $\Ndim_RM=d$, then by \cite[Theorem
4.8]{CN1}, $H^{\fa}_i(M)=0$ for all $i>d$. Also, by
\cite[Proposition 4.2]{Ri}, if $M$ is either finitely generated or
Artinian, then $H^{\fa}_i(M)=0$ for all $i>\sup \{\dim R/\fp|\fp\in
\Cosupp_RM\}$.

When $R$ is complete local, we show that $\Mag_RM\leq \Ndim_RM$ with
the equality if $M\neq 0$ and it is either $N$-critical or semi
discrete linearly compact. So, for having a sharper upper bound for
vanishing of the local homology modules $H^{\fa}_i(M)$, $\Mag_RM$
could be a better candidate. In fact, we conjecture that
$H^{\fa}_i(M)=0$ for all $i>\Mag_RM$. Our investigation on this
conjecture is the core of this paper. We show this conjecture in
several cases.

Namely, we prove that if

$\Coass_RM=\mathcal{A}tt_RM$, or

$M$ is finitely generated, Artinian or Matlis reflexive, or

$M$ is linearly compact, or

$R$ is complete local and $M$ has finitely many minimal coassociated 
prime ideals; or

$R$ is complete local with the maximal ideal $\fm$ and $\fm^n M$ is 
minimax for some integer $n\geq 0$,\\
then $H^\fa_i(M)=0$ for all $i>\Mag_RM$.

Z\"{o}schinger \cite{Z2} has conjectured that
any module over a local ring has finitely many minimal coassociated prime 
ideals. Thus by fourth case, over a complete local ring $R$, Z\"{o}schinger's 
conjecture implies our conjecture.

\section{The results}

In what follows, we denote the faithful exact functor
$\Hom_R(-,\underset{\fm\in \Max R}\oplus E(R/\fm))$ by $(-)^\vee $.
Let $M$ be an $R$-module. A prime ideal $\fp$ of $R$ is said to be a
coassociated prime ideal of $M$ if there is an Artinian quotient $L$
of $M$ such that $\fp=(0:_RL)$. The set of all coassociated prime
ideals of $M$ is denoted by $\Coass_RM$. Also, $\mathcal{A}tt_RM$ is
defined by
$$\mathcal{A}tt_RM:=\{\fp\in \Spec R| \ \fp=(0:_RL) \  \ \text{for
some quotient} \   \ L \ \text{of} \  \  M \}.$$ Clearly,
$\Coass_RM\subseteq \mathcal{A}tt_RM$ and the equality holds if
either $R$ or $M$ is Artinian. More generally, if $M$ is representable, 
then it is easy to check that $\Coass_RM=\mathcal{A}tt_RM$. If $0\lo 
M\lo N\lo L\lo 0$ is
an exact sequence of $R$-modules and $R$-homomorphisms, then it is easy 
to check that
$$\Coass_RL\subseteq \Coass_RN\subseteq \Coass_RL\cup \Coass_RM$$ and
$$\mathcal{A}tt_RL\subseteq \mathcal{A}tt_RN\subseteq \mathcal{A}tt_RL
\cup \mathcal{A}tt_RM.$$ Also, if $R$ is local, then one can see that 
$\Coass_RM=\Ass_RM^{\vee}$.

\begin{lemma} Let $M$ be an $R$-module. Then $\mathcal{A}tt_RM=
\mathcal{A}tt_RM^{\vee\vee}$.
\end{lemma}

\begin{prf} Let $\fp$ be a prime ideal of $R$ and $X,Y$ two $R$-modules.
There are natural isomorphisms $(X/\fp X)^{\vee}\cong
(0:_{X^{\vee}}\fp)$ and $(0:_Y\fp)^{\vee}\cong Y^{\vee}/\fp
Y^{\vee}$. Hence $(X/\fp X)^{\vee\vee}\cong X^{\vee\vee}/\fp
X^{\vee\vee}$. Since $(-)^{\vee}$ is a faithfully exact functor, we
deduce that
$$\Ann_R(X/\fp X)=\Ann_R((X/\fp
X)^{\vee\vee})=\Ann_R(X^{\vee\vee}/\fp X^{\vee\vee}).$$ On the other
hand, one can easily check that $\fp\in \mathcal{A}tt_RX$ if and
only if $\fp=\Ann_R(X/\fp X)$. This yields that
$\mathcal{A}tt_RM=\mathcal{A}tt_RM^{\vee\vee}$, as required.
\end{prf}

Let $\fa$ be an ideal of $R$ and $\mathcal{C}_0(R)$ denote the
category of $R$-modules and $R$-homomorphisms. It is known that the
$\fa$-adic completion functor
$$\Lambda_{\fa}(-):=\underset{n}{\vpl}(R/\fa^n\otimes_R-):
\mathcal{C}_0(R)\lo \mathcal{C}_0(R)$$ is not right exact in
general. For any integer $i$, the $i$-th local homology functor with
respect to $\fa$ is defined as $i$-th left derived functor of
$\Lambda_{\fa}(-)$. For an $R$-module $M$, set $\cd_{\fa}M:=\sup
\{i|H_{\fa}^i(M)\neq 0\}$. By \cite[Corollary 3.2]{GM},
$H^{\fa}_i(M)=0$ for all $i>\cd_{\fa}R$.

\begin{lemma} Let $\fa$ be an ideal of $R$ and $M$ an $R$-module.
Then $H^\fa_i(M)=0$ for all $$i>\sup \{\dim R/\fp|\fp\in 
\mathcal{A}tt_RM\}.$$
\end{lemma}

\begin{prf} For any $R$-module $N$, let $d_N:=\sup \{\dim R/\fp|\fp\in 
\mathcal{A}tt_RN\}$. If $d_N\geq \cd_{\fa}R$, then \cite[Corollary 3.2]{GM} 
implies that $H^\fa_i(N)=0$ for all $i>d_N$. Hence, it is enough to show that
for any $R$-module $N$ with $d_N<\cd_{\fa}R$, one has $H^\fa_i(N)=0$ for all
$d_N<i\leq \cd_{\fa}R+1$. We do decreasing induction on $i$. Clearly, 
the claim holds for $i=\cd_{\fa}R+1$. Now, assume that $d_N<i<\cd_{\fa}R+1$ 
and that the claim holds for $i+1$. We have to show that $H^\fa_i(N)=0$. We
have an exact sequence $$0\lo N\lo N^{\vee\vee} \lo C\lo 0,$$ which
yields the long exact sequence $$\cdots \lo H^\fa_{i+1}(N^{\vee\vee})\lo 
H^\fa_{i+1}(C)\lo H^\fa_i(N)\lo H^\fa_{i}(N^{\vee\vee})\lo \cdots  .$$ 
By \cite[Lemma 3.7]{GM}, one has $H^{\fa}_j(N^{\vee\vee})\cong H^j_{\fa}
(N^{\vee})^{\vee}$ for all $j\geq 0$. It is easy to see that 
$\Ass_RN^{\vee}\subseteq \mathcal{A}tt_RN^{\vee\vee}$, and so by Lemma 2.1,
$\dim_RN^{\vee}\leq d_N$. Thus by Grothendieck's Vanishing Theorem,
one has $$H^\fa_{i+1}(N^{\vee\vee})=0=H^\fa_{i}(N^{\vee\vee}),$$ and
so $H^\fa_{i+1}(C)\cong H^\fa_i(N)$. From the above short exact
sequence and Lemma 2.1, one has $\mathcal{A}tt_RC\subseteq
\mathcal{A}tt_RN$. Hence, $d_C\leq d_N<i+1,$ and so by induction
hypothesis, $$H^\fa_i(N)\cong H^\fa_{i+1}(C)=0.$$ Thus the claim
follows by induction.
\end{prf}

Next, we recall the definitions of $\Ndim_RM$ and $\Mag_RM$.

\begin{definition} Let $M$ be an $R$-module.
\begin{enumerate}
\item[i)] (See \cite{Ro}) The Noetherian dimension of $M$ is defined
inductively as follows: when $M=0$, put $\Ndim_RM=-1$. Then by
induction, for an integer $d\geq 0$, we put $\Ndim_RM=d$ if
$\Ndim_RM<d$ is false and for every ascending sequence $M_0\subseteq
M_1\subseteq\ldots$ of submodules of $M$, there exists $n_0$ such
that $\Ndim_RM_{n+1}/M_n<d$ for all $n>n_0$.
\item[ii)] (See \cite{Y}) The  magnitude of $M$ is defined by $\Mag_RM:=
\sup \{\dim R/\fp|\fp\in \Coass_RM\}$. If $M=0$, we put
$\Mag_RM=-\infty$.
\item[iii)] (See \cite{Ri}) The co-localization of $M$ at a prime ideal
$\fp$ of $R$ is defined by
$$^{\fp}M:=\Hom_{R_{\fp}}((M^\vee)_{\fp},E_{R_{\fp}}({R_{\fp}}/\fp
{R_{\fp}})).$$ Then $\Cosupp_RM$ is defined by $\Cosupp_RM:=\{\fp\in
\Spec R|^{\fp}M\neq 0\}$.
\item[iv)] (See \cite{Cha}) $M$ is said to be $N$-critical if $\Ndim_RN
<\Ndim_RM$ for all proper submodules $N$ of $M$.
\end{enumerate}
\end{definition}

From the definition, it becomes clear that $\Ndim_RM=0$ if and only if 
$M$ is a non-zero Noetherian $R$-module. If $0\lo X\lo Y\lo Z\lo 0$ is an 
exact sequence of $R$-modules and $R$-homomorphisms, then 
\cite[Proposition 5]{Cha} yields that $\Ndim_RY=\max\{\Ndim_RX,
\Ndim_RZ\}$. Also, it is easy to verify that $\Mag_RY=\max\{\Mag_RX,
\Mag_RZ\}$.

Next, we compare $\Ndim_RM$ and $\Mag_RM$. Recall that an $R$-module $M$ 
is said to be {\it Matlis reflexive} if the natural homomorphism $M\lo 
M^{\vee\vee}$ is an isomorphism.

\begin{lemma} Let $M$ be an $R$-module.
\begin{enumerate}
\item[i)] Suppose $R$ is complete local. Then $\Mag_RM\leq \Ndim_RM$ and
equality holds if $M\neq 0$ and it is either $N$-critical or Matlis
reflexive.
\item[ii)] $\Mag_RM\leq \sup \{\dim R/\fp|\fp\in \Cosupp_RM\}$ and 
equality holds if $R$ is local.
\end{enumerate}
\end{lemma}

\begin{prf} i) Let $\fp\in\Coass_RM$. Then there is an Artinian $\fp$-secondary 
quotient $M/N$ of $M$ such that $\fp=\Ann_RM/N$. By \cite[Proposition 5]{Cha},
we have $\Ndim_RM=\max \{\Ndim_RN, \Ndim_RM/N\}$. On the other hand, for any 
Artinian $R$-module $A$, \cite[Theorem 2.10]{Y} asserts that $\Mag_RA=\Ndim_RA$. 
(Note that the argument of \cite[Theorem 2.10]{Y} is not correct without
the completeness assumption on $R$.) Since $\fp$ is the only coassociated prime 
ideal of $M/N$, it turns out that $$\dim R/\fp=\Mag_RM/N=\Ndim_RM/N\leq \Ndim_RM.$$ 
Thus $$\Mag_RM=\sup \{\dim R/\fp|\fp\in \Coass_RM\}\leq\Ndim_RM.$$

Now, assume that $M\neq 0$. Let $M$ be $N$-critical and let $\fp$ and $N$ be as 
above. Then the $\Ndim_RN<\Ndim_RM$, and so
$$\begin{array}{llll} \Ndim_RM&=\max\{\Ndim_RN,\Ndim_RM/N\}\\
&=\Ndim_RM/N\\
&=\Mag_RM/N\\
&\leq \Mag_RM.
\end{array}
$$
Next, let $M$ be Matlis reflexive. Then, by \cite[Theorem
12]{BEG}, there is a finitely generated submodule $N$ of $M$ such
that $M/N$ is Artinian. Clearly, we may suppose that $M\neq 0$.
Then $$\Ndim_RM=\max\{0,\Ndim_RM/N\}=\max\{0,
\Mag_RM/N\}\leq \Mag_RM.$$

ii) It suffices to show that $\Coass_RM\subseteq \Cosupp_RM$. Let
$\fp\in \Coass_RM$. There is an Artinian quotient $L$ of $M$ such
that $\fp=\Ann_RL$. Set $E:=\oplus_{\fm\in \Max R}E(R/\fm)$. Since
$L$ is Artinian, we may assume that $L\subseteq E^r$ for an integer
$r\geq 0$. Let $f:M\lo E^r$ denote the composition of the natural
epimorphism $M\lo L$ and the natural monomorphism $L\lo E^r$. Then
$\fp=(0:_Rf)(:=\{a\in R|af:M\lo E^r \   \  \text{is the zero
homomorphism} \})$, and so $$\fp\in
\Ass_R(\Hom_R(M,E^r))=\Ass_RM^\vee\subseteq \Supp_RM^\vee.$$ Hence,
$(M^\vee)_{\fp}\neq 0$, and so $^{\fp}M\neq 0$. This means that
$\fp\in \Cosupp_RM$.

Now, assume that $R$ is local. Since $\Cosupp_RM=\Supp_RM^\vee$,
\cite[Lemma 2.2 a)]{Y} implies that $\Mag_RM=\sup \{\dim R/\fp|\fp\in
\Cosupp_RM\}$.
\end{prf}

Recall that an $R$-module $M$ is said to be {\it semi-Artinian} if every proper
submodule of $M$ contains a minimal submodule; see e.g. \cite{Ru}.

\begin{example} Let $M$ be an $R$-module.
\begin{enumerate}
\item[i)] Suppose that $M$ is semi-Artinian with finitely many associated prime ideals.
Then by \cite[Bemerkung after Satz 2.9]{Z2},  one has $\Coass_RM=\mathcal{A}tt_RM$. But,
in general, the containment $\Coass_RM\subseteq \mathcal{A}tt_RM$ may be strict. To this
end, let $\fp$ be a non-maximal prime ideal of $R$ and set $M:=R/\fp$. Then $\Coass_RM=\V(\fp)
\cap \Max R$ and $\mathcal{A}tt_RM=\V(\fp)$, and so  $\Coass_RM\varsubsetneq \mathcal{A}tt_RM$.
\item[ii)] The inequality in Lemma 2.4 i) may be strict. To see this, let
$(R,\fm)$ be a local ring and $M:=\oplus_{i\in \mathbb{N}} R/\fm$.
Then $\Coass_RM=\{\fm\}$, and so $\Mag_RM=0$. But, as $M$ is not a
Noetherian $R$-module, we have $\Ndim_RM>0$.
\item[iii)] In Lemma 2.4 i), the completeness assumption on $R$ can not be
skipped. To see this, let $(R,\fm)$ be a two-dimensional local
domain such that $\hat {R}$ possesses a one-dimensional embedded
prime ideal $\fq$; see \cite[12, Appendix, Example 2]{N}. Let
$A:=(0:_{E(R/\fm)}\fq)$. Then $\Mag_RA=2$ and $$\Ndim_RA=\Ndim_{\hat{R}}A=\Mag_{\hat{R}}A=1.$$
\end{enumerate}
\end{example}

\begin{lemma} Let $(R,\fm)$ be a complete local ring, $\fa$ an ideal of
$R$ and $M$ an $R$-module. Then $H_i^\fa(M)=0$ for all $i>\Mag_RM^{\vee\vee}$.
\end{lemma}

\begin{prf} The proof is similar to the proof of Lemma 2.2.
We use decreasing induction on $i$. For $i\geq \dim R+1$, the
claim holds by \cite[Corollary 3.2]{GM}. Note that
$\Mag_RM^{\vee\vee}\leq \dim R$. Now, assume that
$\Mag_RM^{\vee\vee}<i<\dim R+1$ and that the claim holds for $i+1$.
We have an exact sequence
$$0\lo M\lo M^{\vee\vee} \lo C\lo 0,$$ which yields the long exact
sequence
$$\cdots \lo H^\fa_{i+1}(M^{\vee\vee})\lo H^\fa_{i+1}(C)\lo H^\fa_i(M)\lo H^\fa_{i}
(M^{\vee\vee})\lo \cdots  \ \ .$$ By \cite[Lemma 3.7]{GM}, one has
$H^{\fa}_j(M^{\vee\vee})\cong {H^j_{\fa}(M^{\vee})}^{\vee}$ for all
$j\geq 0$. Since $$i>\Mag_RM^{\vee\vee}=\dim_R M^{\vee\vee\vee}\geq
\dim_RM^{\vee},$$ Grothendieck's Vanishing Theorem implies that
$$H^\fa_{i+1}(M^{\vee\vee})=0=H^\fa_{i}(M^{\vee\vee}).$$ Hence
$H^\fa_{i+1}(C)\cong H^\fa_i(M)$. Also, from the above short exact
sequence, we deduce that
$$\Mag_R(M^{\vee\vee\vee\vee})=\max \{\Mag_RM^{\vee\vee},\Mag_RC^{\vee\vee}\}.$$
On the other hand, \cite[Lemma 2.9 and Folgerung 2.10]{Z1} yields
that $\Coass_R(M^{\vee\vee\vee\vee})=\Coass_R(M^{\vee\vee})$, and so
$$\Mag_RC^{\vee\vee}\leq
\Mag_R(M^{\vee\vee\vee\vee})=\Mag_RM^{\vee\vee}.$$ Now, since
$\Mag_RC^{\vee\vee}<i+1$, by induction hypothesis, it turns out that
$$H^\fa_i(M)\cong H^\fa_{i+1}(C)=0.$$ Thus, the claim follows by
induction.
\end{prf}

\begin{lemma} Let $(R,\fm)$ be a complete local ring and $M$ an
$R$-module. Assume that $M$ has finitely many minimal coassociated
prime ideals. Then $\Mag_RM^{\vee\vee}=\Mag_RM$.
\end{lemma}

\begin{prf} Let $\{\fp_1,\fp_2,\dots ,\fp_n\}$ be the set of all elements of
$\Coass_RM$ which are minimal with respect to inclusion in
$\Coass_RM$. By \cite[Satz 2.6]{Z1}, for any prime ideal
$\tilde{p}$ of $R$, $\tilde{p}\in \Coass_RM^{\vee\vee}$ if and only
if $\tilde{p}=\bigcap_{\fp\in \Lambda}\fp$ for some subset $\Lambda$
of $\Coass_RM$. In particular, one has $\Coass_RM\subseteq
\Coass_RM^{\vee\vee}$, and so $\Mag_RM\leq \Mag_RM^{\vee\vee}$. Also, it follows
that any prime ideal $\tilde{p}\in \Coass_RM^{\vee\vee}$ contains $\underset{i=1}
{\overset{n} \bigcap}\fp_i$.

Let $\fq\in \Coass_RM^{\vee\vee}$ be such that $\dim
R/\fq=\Mag_RM^{\vee\vee}$. Then $\fq$ is minimal in
$\Coass_RM^{\vee\vee}$. Since $\fq\supseteq
\underset{i=1}{\overset{n} \bigcap}\fp_i$, there is $1\leq j\leq n$
such that $\fq\supseteq \fp_j$. But $\fp_j\in\Coass_RM^{\vee\vee}$,
and so $\fq=\fp_j$. Thus $\fq\in \Coass_RM$, and so
$\Mag_RM^{\vee\vee}\leq \Mag_RM$.
\end{prf}

At this point, we are ready to present our main result. But, we first recall some
definitions which are needed in its statement.

We begin by recalling the definition of linearly compact modules from \cite{Mac}. Let $M$
be a topological $R$-module. Then $M$ is said to be {\it linearly topologized} if $M$ has a base
$\mathcal{M}$ consisting of submodules for the neighborhoods of its zero element. A Hausdorff
linearly topologized $R$-module $M$ is said to be {\it linearly compact} if for any family
$\mathcal{F}$ of cosets of closed submodules of $M$ which has the finite intersection property,
the intersection of all cosets in $\mathcal{F}$ is non-empty. A Hausdorff linearly
topologized $R$-module $M$ is called {\it semi-discrete} if every submodule of $M$ is closed. The class
of semi-discrete linearly compact modules is very large, it contains many important classes of
modules such as the class of Artinian modules, or the class of finitely generated modules over a
complete local ring.

An $R$-module $M$ is called {\it minimax} if it has a finitely generated submodule $N$ such that
$M/N$ is Artinian. By \cite[Lemma 1.1]{Z1}, over a complete local ring $R$, an $R$-module $M$ is
minimax if and only if $M$ is semi-discrete linearly compact and if and only if $M$ is Matlis
reflexive.

\begin{theorem} Let $\fa$ be an ideal of $R$ and $M$ an $R$-module. Assume that either
\begin{enumerate}
\item[i)] $\Coass_RM=\mathcal{A}tt_RM$,
\item[ii)] $M$ is $N$-critical,
\item[iii)] $M$ is finitely generated, Artinian or Matlis reflexive,
\item[iv)] $M$ is linearly compact,
\item[v)] $R$ is complete local and $M$ has finitely many minimal coassociated prime ideals; or
\item[vi)] $R$ is complete local with the maximal ideal $\fm$ and $\fm^n M$ is minimax for some
integer $n\geq 0$.
\end{enumerate}
Then $H^\fa_i(M)=0$ for all $i>\Mag_RM$.
\end{theorem}

\begin{prf} i) follows by Lemma 2.2.

ii) By \cite[Proposition 2]{Cha}, $M$ is a secondary module. This implies that
$\Coass_RM=\mathcal{A}tt_RM$, and so ii) follows by i).

iii) For a finitely generated $R$-module $N$, one has $\Mag_RN\leq 0$ and
\cite[Proposition 3.2]{Si1} yields that $H^\fa_i(N)=0$ for all
$i>0$. When $M$ is Artinian, the claim follows by i).  Now, assume
that $M$ is Matlis reflexive. We may and do assume that $M$ is nonzero.
By \cite[Theorem 12]{BEG}, there is a finitely generated submodule $N$
of $M$ such that $M/N$ is Artinian. Then
$$\Mag_RM=\max\{\Mag_RN,\Mag_RM/N\}=\max\{0,\Mag_RM/N\}.$$
Thus, the claim follows by the following long exact sequence
$$\cdots \lo H^\fa_{i+1}(M/N)\lo H^\fa_i(N)\lo H^\fa_i(M)\lo H^\fa_{i}
(M/N)\lo \cdots  \ \ .$$

iv) Let $\mathcal{M}$ be a base consisting of submodules for the neighborhoods of
the zero element of $M$.  Then by \cite[3.11]{Mac}, $M\cong \underset{U\in
\mathcal{M}}{\vpl}M/U$. By \cite[Proposition 3.4]{CN1}, it turns out that
$H^\fa_i(M)\cong \underset{U\in \mathcal{M}}{\vpl}H^\fa_i(M/U)$. Note that
for each $U$, $M/U$ is a semi discrete linearly compact
$R$-module with $\Mag_RM/U\leq \Mag_RM$. So, we only need to prove the claim for
the case $M$ is a semi discrete linearly compact $R$-module. Let $M$ be a  semi discrete
linearly compact $R$-module. Then \cite{Z3} implies that $M$ is Matlis reflexive. Thus iii)
completes the proof of this part.

v) follows by Lemmas 2.6 and 2.7.

vi) By \cite[Theorem 3.3]{Ru} and \cite[Lemma 2.2 d)]{Z3}, there is a submodule $N$ of
$M$ such that $\Coass_RN\subseteq \{\fm\}$ and the quotient module $M/N$ is Artinian.
Hence $\Coass_RM$ is finite, and so v) yields the conclusion.
\end{prf}

Let $(R,\fm)$ be a local ring. Z\"{o}schinger \cite{Z2} conjectured that any
$R$-module $M$ has finitely many minimal coassociated prime ideals. Clearly, this
is equivalent to say that for any $R$-module $M$, $\Supp_RM^{\vee}$ is a
Zariski-closed subset of $\Spec R$.

\begin{remark} Let $M$ be an $R$-module and in the first three parts, suppose that $R$ is local.
\begin{enumerate}
\item[i)] By \cite[Folgerung 1.5]{Z2} if $\Coass_RM$ is countable, then $M$ has finitely many
minimal coassociated prime ideals. Hence, if $R$ is countable, or $\Spec R$ is finite (e.g.
$\dim R\leq 1$), or $M$ is Matlis reflexive; or $M$ is representable (e.g. injective), then $M$
has finitely many minimal coassociated prime ideals.
\item[ii)] Let $\{M_{\lambda}\}_{\lambda\in \Lambda}$ be a family of $R$-modules such that for
each $\lambda$, $M_{\lambda}$ has finitely many minimal coassociated prime ideals. Then both of $\bigoplus_{\lambda\in \Lambda}M_{\lambda}$ and $\prod_{\lambda\in \Lambda}M_{\lambda}$ have
finitely many minimal coassociated prime ideals; see \cite[Satz 2.6]{Z2} and \cite[Satz 2.8 b]{Z1}.
\item[iii)] If $\Coass_RM=\mathcal{A}tt_RM$, then by \cite[Lemma 3.1]{Z2} $M$ has finitely many
minimal coassociated prime ideals. Hence for any infinite index set $\Lambda$ and any $R$-module
$X$, the $R$-modules $X^{(\Lambda)}$ and $X^{\Lambda}$ have finitely many minimal coassociated
prime ideals; see \cite[Bemerkung after Satz 2.4]{Z2}.
\item[iv)] In view of Lemma 2.4 ii),  clearly Theorem 2.8 iii) extends \cite[Proposition 4.2]{Ri}.
\end{enumerate}
\end{remark}

\begin{acknowledgement} We would like to thank Professor Helmut  Z\"{o}schinger for his comments
on coassociated prime ideals. Part of this research was done during the second author's visit to
the Department of Mathematics at the University of Nebraska-Lincoln. He thanks this department for
its kind hospitality.
\end{acknowledgement}


\end{document}